%% file: ims91-19b.tex
\documentstyle[fullpage,12pt]{article}

\input{psfig}
\input{imsmark}
\begin{document}
\thispagestyle{empty}
\vspace*{1in}

\def\IMSmarkvadjust{-.25in}
\SBIMSMark{1991/19b}{June 1991}{}

\centerline {\bf On the quasisymmetrical classification of infinitely
renormalizable maps}

\vskip5pt
\centerline {\sc II. Remarks on maps with a bounded type topology}

\vskip10pt
\centerline{by Yunping Jiang }

\vskip5pt
\centerline{June, 1991}

\vskip15pt
\noindent {\bf \S 0 Introduction}

\vskip5pt
This note is a remark to the paper \cite{j1}. 
The aim is to
show that the techniques in \cite{j1} can also be used to understand the
quasisymmetrical classification of infinitely renormalizable maps of 
bounded type.   
We will use the same terms and notations as those in \cite{j1} without further
notices. The result we will prove is the following theorem.

\vskip5pt
{\sc Theorem 1.} {\em 
Suppose $f$ and $g$ in ${\cal U}$ are two infinitely
renormalizable maps of bounded type and topologically conjugate. 
Moreover, suppose $H$ is
the homeomorphism between $f$ and $g$. Then $H$ is quasisymmetric.}

\vskip5pt
Since the techniques as well as ideas of the proof are similar to those in
\cite{j1}, we outline the proof in the next section. The reader may
refer to \cite{j1} and \cite{j2} for more details.

\vskip15pt
\noindent {\bf \S 1 The outline of the proof of Theorem 1}

\vskip5pt
We outline the proof of Theorem 1 by several lemmas.

\vskip5pt
Suppose $f=h\circ Q_{t}$, for some $t>1$, in ${\cal U}$ is an
infinitely renormalizable map of bounded type. We note that
$Q_{t}(x) = -|x|^{t}$. Let $f_{0}=f$. And inductively, 
let $f_{k} = \alpha_{k}\circ f_{k-1}^{\circ
n_{k}}\circ \alpha_{k}^{-1}$ be the renormalization ${\cal R}(f_{k-1})$ 
of $f_{k-1}$ where
$\alpha_{k}$ is the linear rescale from $J_{k}$ to $[-1,1]$ and $n_{k}$ is
the return time for any integer
$k\geq 1$ (see \cite{j1}). We call $J_{k}$ a restricted interval.

\vskip5pt 
Let $I_{0}$ be the interval $[-1,1]$ and $I_{k}$ be the preimage of
$[-1, 1]$ under $\alpha_{1}\circ \cdots \circ \alpha_{k}$ for $k\geq 1$. 
We note that the set $\{ I_{k} \}_{k=0}^{\infty}$ forms a sequence of
nested intervals.  Moreover, one of the endpoints of $I_{k}$, say
$p_{k}$, is a periodic
point of period $m_{k} =n_{1}\cdots n_{k}$ of $f$ and the orbit $O(p_{k})$ of $p_{k}$ under
$f$ stays outside of the interior of $I_{k}$ (see Figure 1).

\vskip5pt
\centerline{\psfig{figure=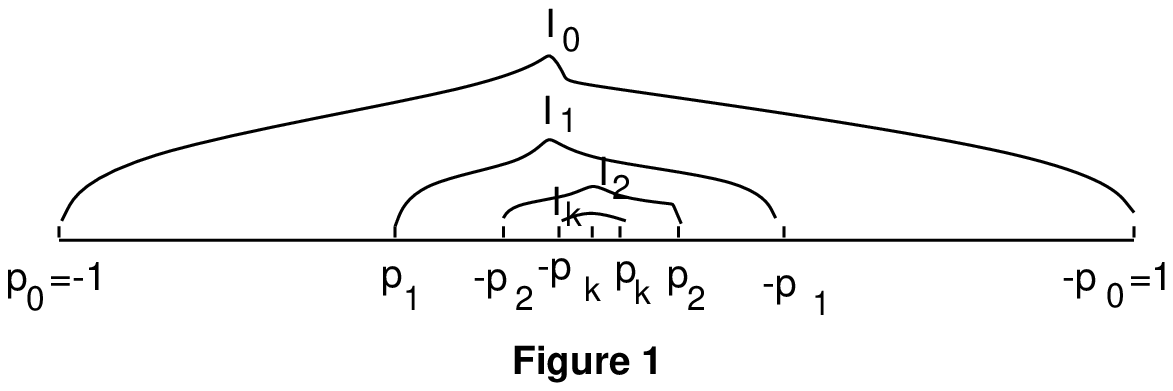}}

\vskip5pt
Suppose $L_{k}$ is the image of $I_{k}$ under $f^{\circ m_{k}}$ and
$T_{k}$ is the interval bounded by the points $p_{k}$ and $p_{k+1}$. Let 
$M_{k}$ be the complement of $T_{k}$ in $L_{k}$. Then $M_{k}$
is the interval bounded by $p_{k+1}$ and $c_{m_{k}}$, where 
$c_{m_{k}}=f^{\circ m_{k}}(0)$ (See Figure 2). 
  
\vskip5pt
\centerline{\psfig{figure=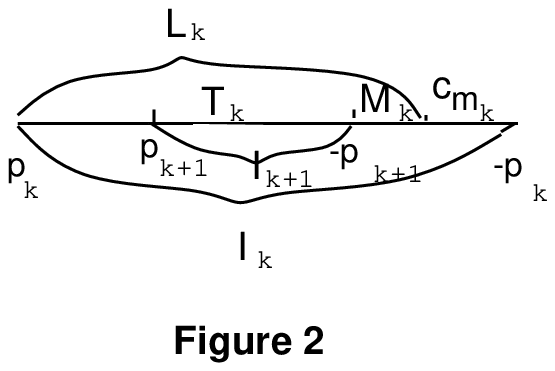}}

\vskip5pt 
{\sc Lemma 1.} {\em There is a constant $C_{1}=C_{1}(f)>0$ such that
\[ C_{1}^{-1} \leq |M_{k}|/|I_{k}|\leq C_{1}. \]
for all the integers $k \geq 0$ .} 

\vskip5pt
{\sl Proof.} This lemma is actually proved in \cite{s1} by using the techniques
such as the smallest
interval and shuffle permutation on the intervals. 

\vskip5pt
{\sc Lemma 2.} {\em There is a constant $C_{2}=C_{2}(f)>0$ such that
\[ C_{2}^{-1} \leq |I_{k}|/|I_{k-1}|\leq C_{2}\]
for all the integers $k\geq 0$.}
 
\vskip5pt
We first prove a more general result, as that in \cite{j1}, as follows.
Let ${\cal K}= {\cal K}(t,N,K)$, for fixed numbers $t>1$, $N\geq 2$ and
$K>0$, be the subspace of renormalizable 
maps $f=h\circ Q_{t}$ in ${\cal U}$ such that $|\Big( N(h)\Big) (x)| \leq K$
for all $x$ in $[-1,0]$ 
and all the return times $n_{k}$ of ${\cal R}^{\circ k}(f)$ are less than or
equal to $N$.

\vskip5pt
{\sc Lemma 3.} {\em There is a constant $C_{3}=C_{3}(t,N,K)>0$ such that
\[ C_{3}^{-1} \leq f(0)=c_{1}(f) \leq C_{3}\] 
for all $f$ in ${\cal K}$.}

\vskip5pt
{\sl Proof.} 
The proof of this lemma is similar to the proof of Lemma 3 in \cite{j1} but needs
little more work to solve a little more complicated equation.

\vskip5pt
Remember that $f_{k}$ is the $k^{th}$-renormalization of $f=f_{0}$. Let $f_{k}= h_{k}\circ Q_{t}$. 
We note that the graph of $f_{k}$
is the rescale of the graph of the restriction of $f^{\circ m_{k}}$ to
$I_{k}$. 

\vskip5pt
{\sc Lemma 4}. {\em There is a constant $C_{4}=C_{4}(f)>0$ such that  
\[ |\Big( N(h_{k})\Big) (x)| \leq C_{4} \]
for all $x$ in $[-1, 0]$ and all the integer $k\geq 0$.}

\vskip5pt
{\sl Proof.} It is the a prior bound proved in \cite{s1}.

\vskip5pt
{\sl Proof of Lemma 2}. It is now a direct corollary of 
Lemma 1, Lemma 3 and Lemma 4 for $K= C_{4}$ and $N=\max_{0\leq k< \infty}\{
n_{k} \}$.

\vskip5pt
The set of the nested intervals $\{ I_{0}$, $I_{1}$, $\cdots $, $I_{k}$,
$\cdots \}$ gives a partition of $[-1,1]$ as follows.
Let $p_{k}$ be one of the endpoints of $I_{k}$ and $O_{k,f}(p_{k})$ be
the intersection of $I_{k-1}$ and the orbit $O_{k, f}(p_{k})$ of 
$p_{k}$ under $f$ for
$k\geq 1$.
Suppose $M_{k-1, 1}$, $\cdots$, $M_{k-1, n_{k}+1}$ are the connected 
components of
$I_{k-1}\setminus (O_{k, f}(p_{k})\cup I_{k})$ for $k\geq 1$.  
Then the set $\eta_{0}= \{ M_{k, i} \}$ for $i=1$, $\cdots $, $n_{k}+1$ and
$k=1$, $2$, $\cdots $ forms a partition of $[-1,1]$ (see Figure 3).

\vskip5pt
\centerline{\psfig{figure=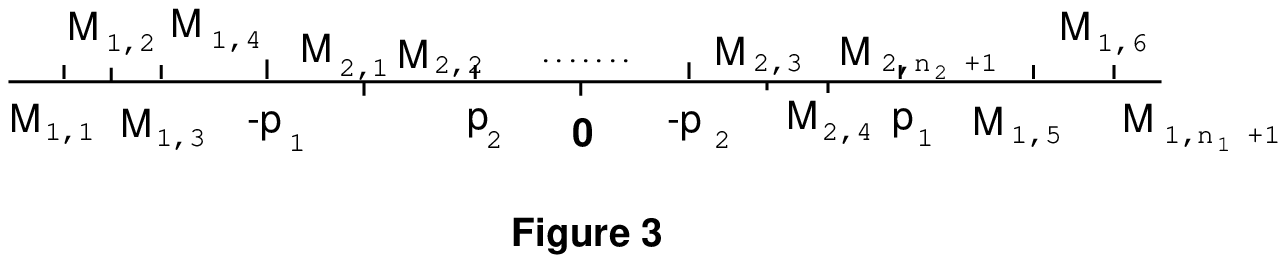}}

\vskip5pt
Now we are going to define a Markov map $F$ induced by $f$. Let $F$ be a
function of $[-1,1]$ defined by
\[ F(x) =\left\{ \begin{array}{ll}
			 f(x) ,& x\in M_{1,1} \cup M_{1,2} \cup \cdots \cup M_{1, n_{1}+1},\\
			 f^{\circ n_{1}}(x) ,& x\in M_{2,1}\cup \cdots \cup M_{2, n_{2}+1},\\
                         \vdots &  \\
                         f^{\circ n_{1}n_{2}\cdots n_{k}}(x), & x\in M_{k,
1}\cup \cdots \cup
M_{k, n_{k}+1},\\
                         \vdots & 
			 \end{array}
		\right.\]

\vskip5pt
It is clearly that $F$ is a Markov map in the sense that the image of every
$M_{k, i}$ under $F$ is the union of some intervals in $\eta_{0}$ (Figure 4).

\vskip5pt
Let $g_{k, i}=(F|M_{k, i})^{-1}$ for $k=1$, $\cdots $, and $i=1$, $\cdots
n_{k}+1$ be the inverse
branches of $F$ with respect to the Markov partition $\eta_{0}$. Suppose $w=i_{0}i_{1}\cdots i_{l-1} $ is a finite sequence of
the set ${\cal I}= \{ (k,i)$, $ k=1$, $\cdots$ and $i=1$, $\cdots $, $n_{k}+1 \}$. 
We say it is admissible if the range $M_{i_{s}}$ of
$g_{i_{s}}$ is contained in the domain $F_{i_{s-1}}(J_{i_{s-1}})$ of
$g_{i_{s-1}}$ for $s=1$, $\cdots $, $l-1$.
For an admissible sequence $w = i_{0}i_{1}\cdots i_{l-1}$, we can define
$g_{w}=g_{i_{0}}\circ g_{i_{1}}\circ  \cdots \circ g_{i_{l-1}}$.  We use
$D(g_{w})$ to denote the domain of $g_{w}$ and use $|D(g_{w})|$ to denote
the length of the interval $D(g_{w})$. 

\vskip5pt
\centerline{\psfig{figure=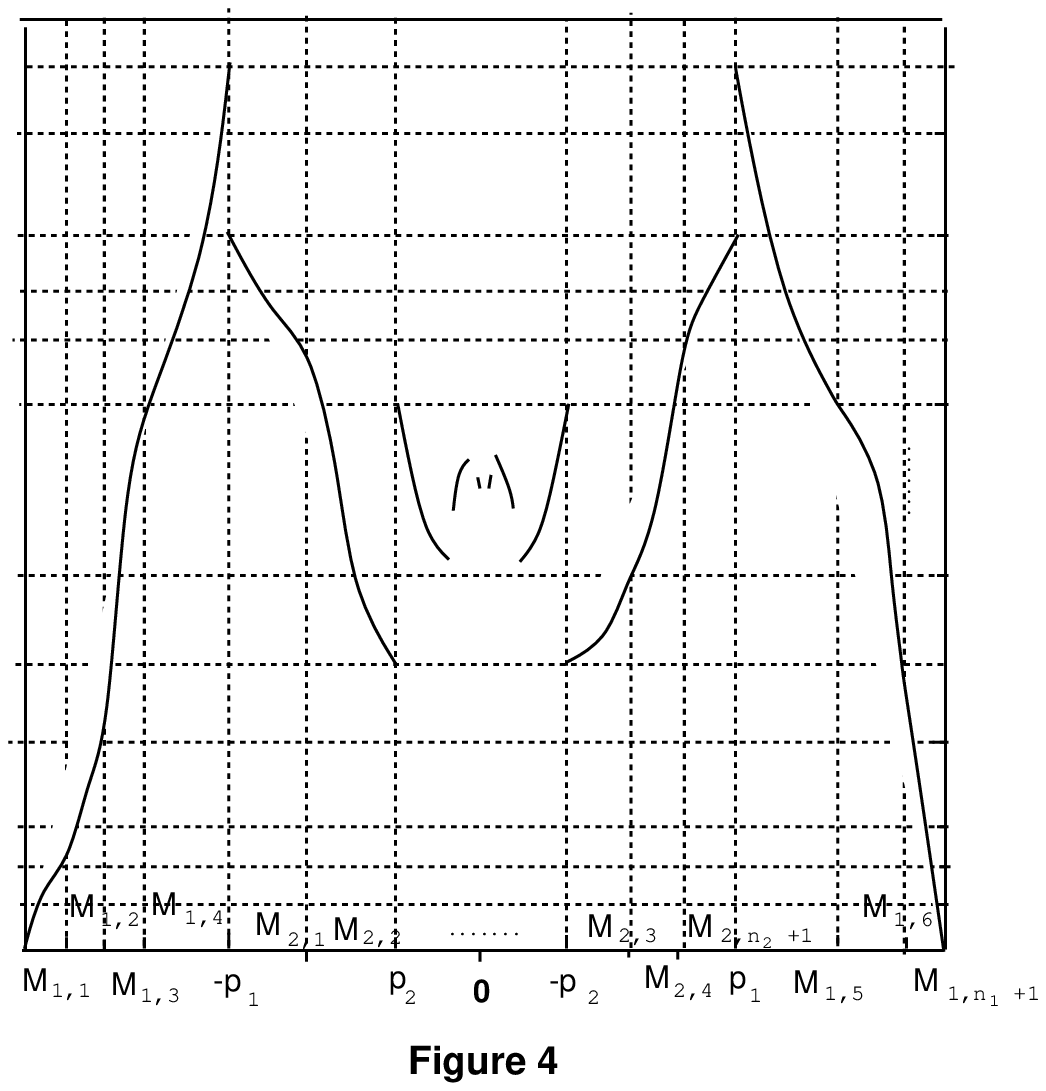}}

\vskip5pt
{\sc Definition 1.} {\em We say the induced Markov map $F$ has bounded distortion
property if there is a constant $C_{5}=C_{5}(f)>0$ such that 

\vskip5pt
(a) $C_{5}^{-1}\leq |M_{k, i}|/ |M_{k, i+1}|\leq C_{5}$ for $k=1$, $2$,
$\cdots $ and $i=1$, $\cdots
$, $n_{k}$,

(b) $C_{5}^{-1} \leq |M_{k, i}|/|I_{k}| \leq C_{5}$ for $k=1$, $2$, $\cdots$
and $i=1$,$\cdots$,
$n_{k}+1$,  and 

\vskip5pt
(b) $| \Big( N(g_{w})\Big) (x)|\leq C_{5}/|D(g_{w})|$ for
all $x$ in $D(g_{w})$ and all admissible $w$.} 

\vskip5pt
The reason we give this definition is the following lemma as that in \cite{j1}.

\vskip5pt
{\sc Lemma 5.} {\em Suppose $f$ and $g$ in ${\cal U}$ are two infinitely 
renormalizable
maps of bounded type and $H$ is the conjugacy 
between $f$ and $g$. If both of the induced Markov maps $F$ and $G$ have 
the bounded
distortion property, then $H$ is quasisymmetric.}

\vskip5pt
{\sl Proof.} It can be proved by almost the same arguments as that 
we used in the paper \cite{j2}. For
more details of the proof, the reader may refer to \cite{j2}.    

\vskip5pt
Now the proof of Theorem 1 concentrates on the next lemma.

\vskip5pt
{\sc Lemma 6.} {\em Suppose $f=h \circ Q_{t}$, for some $t>1$, in 
${\cal U}$ is an infinitely renormalizable
map of bounded type
and
$F$ is the Markov map induced by $f$. Then $F$ has 
the bounded distortion property.}

\vskip5pt
{\sl Proof.} Let $I_{k,j}=\Big( f^{\circ m_{k-1}}|I_{k-1}\Big)^{\circ
j}(I_{k})$ for $j=0$, $1$, $\cdots $, $n_{k}$ and $\{ G_{k, i} \}$ are all
the connected components of $I_{k}\setminus \cup_{j=0}^{n_{k}} I_{k, j}$
(Figure 5).

\vskip5pt

Each $M_{k, j}$ is either a single $G_{k, i}$ or $I_{k, j}\cup G_{k, i}$ for
some $j$ and some $i$. By the bounded geometry [3] of $\{ I_{k, j} \}$ and
$\{ G_{k, i} \}$, there is a constant $C_{6} >1$ such that 
all the ratios $|I_{k, j}|/|I_{k, j'}|$, $|G_{k, i}|/|G_{k, i'}|$ and $|G_{k,
i}|/|I_{k, j}|$ are in the interval $[C_{6}^{-1}, C_{6}]$. We note that
$C_{6}$ does not depend on $k$ as well as $i$, $i'$, $j$ and $j'$.
This fact and Lemma 3 imply the condition (a) in Definition 1. 

\vskip5pt
\centerline{\psfig{figure=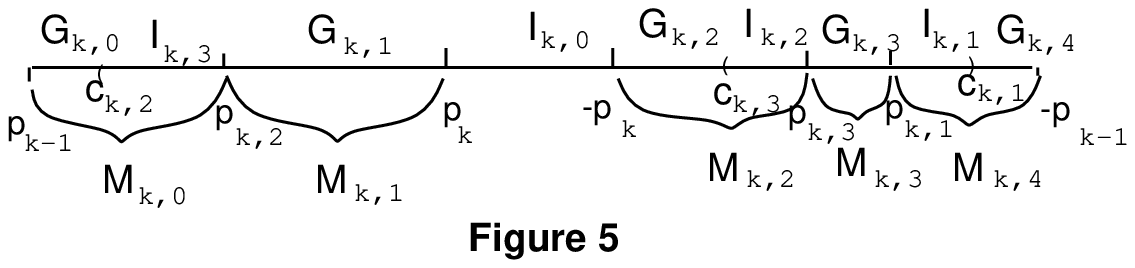}}

\vskip5pt
The condition (b) in Definition 1 is assured by Lemma 2 and the condition (a).
The proof of the condition (c) in Definition 1 is similar to that in \cite{j1}.

\vskip5pt
The arguments in Lemma 1 to Lemma 6 give the proof of Theorem 1.

\vskip20pt
Yunping Jiang

Institute for Mathematical Sciences 

SUNY at Stony Brook

Stony Brook, NY 11794

e-mail: jiang@math.sunysb.edu

\end{document}

%% file: imsmark.tex
\def\IMSmarkvadjust{0 pt}
\def\IMSmarkhadjust{0 pt}
\def\IMSmarkhpadding{0 pt}
\def\IMSpubltext{Published in modified form:}
\def\SBIMSMark#1#2#3{
 \font\SBF=cmss10 at 10 true pt
 \font\SBI=cmssi10 at 10 true pt
 \setbox0=\hbox{\SBF \hbox to \IMSmarkhpadding{\relax}
                Stony Brook IMS Preprint \##1}
 \setbox2=\hbox to \wd0{\hfil \SBI #2}
 \setbox4=\hbox to \wd0{\hfil \SBI #3}
 \setbox6=\hbox to \wd0{\hss
             \vbox{\hsize=\wd0 \parskip=0pt \baselineskip=10 true pt
                   \copy0 \break%
                   \copy2 \break%
                   \copy4 \break}}
 \dimen0=\ht6   \advance\dimen0 by \vsize \advance\dimen0 by 8 true pt
                \advance\dimen0 by -\pagetotal
	        \advance\dimen0 by \IMSmarkvadjust
 \dimen2=\hsize \advance\dimen2 by .25 true in
	        \advance\dimen2 by \IMSmarkhadjust

  \openin2=publishd.tex
  \ifeof2\setbox0=\hbox to 0pt{}
  \else 
     \setbox0=\hbox to 3.1 true in{
                \vbox to \ht6{\hsize=3 true in \parskip=0pt  \noindent  
                {\SBI \IMSpubltext}\hfil\break
                \input publishd.tex 
                \vfill}}
  \fi
  \closein2
  \ht0=0pt \dp0=0pt
 \ht6=0pt \dp6=0pt
 \setbox8=\vbox to \dimen0{\vfill \hbox to \dimen2{\copy0 \hss \copy6}}
 \ht8=0pt \dp8=0pt \wd8=0pt
 \copy8
 \message{*** Stony Brook IMS Preprint #1, #2. #3 ***}
}